\documentclass[11pt]{article}

\usepackage{amsmath, amssymb, amsthm, graphicx}
\usepackage{epstopdf}% To incorporate .eps illustrations using PDFLaTeX, etc.

\usepackage[caption=false]{subfig}% Support for small, `sub' figures and tables

\usepackage[numbers,sort&compress]{natbib}% Citation support using natbib.sty
\bibpunct[, ]{[}{]}{,}{n}{,}{,}% Citation support using natbib.sty
\theoremstyle{plain}% Theorem-like structures provided by amsthm.sty
\newtheorem{theorem}{Theorem}[section]
\newtheorem{lemma}[theorem]{Lemma}

\theoremstyle{definition}

\newtheorem{example}[theorem]{Example}

\theoremstyle{remark}

\begin{document}
% \articletype{ARTICLE TEMPLATE}% Specify the article type or omit as appropriate

\title{Exit times for some nonlinear autoregressive processes}
\author{
G\"oran H\"ogn\"as\thanks{CONTACT G\"oran H\"ogn\"as. Email: ghognas@abo.fi} \;\; and Brita Jung}
\date{}
\maketitle

\begin{center}
Department of Natural Sciences, \AA bo Akademi,\\ FIN-20500 \AA bo, Finland,

\medskip

18.12.2019
\end{center}

\begin{abstract}
By using the large deviation principle, we investigate the expected exit time from the interval $[-1,1]$ of a process $\{X_n^\varepsilon\}$, $n = 0,1,2,\ldots$, of autoregressive type, defined recursively by
\begin{equation}
X_{n+1}^\varepsilon = f(X_n^\varepsilon) + \varepsilon\xi_{n+1}, X_0 = 0.
\end{equation}
Here, $\varepsilon$ is a small positive parameter, $f: R\mapsto R$ is usually a contractive function and $\{\xi_n\}_{n\ge 1}$ is a sequence of i.i.d. random variables. The case when $f$ is linear and $\xi_n$ has a normal distribution has been treated before. In this paper, we extend the results to more general functions $f$, with the main focus on piecewise linear functions.
\end{abstract}

% \begin{keywords}
% Large deviation principle, autoregressive process, exit time
% \end{keywords}

%%%%%%%%%%%%%%%%%%%%%%%%%%%%%%%%%%%%%%%%%%%%%%%%%%%%%%%%%%%%%%%%%%%%%%%%%%%%%%%%%%%%
\section{Introduction}
The original motivation behind our work is to study the time until extinction of a population. A stochastic process that models a population may stay at a certain level (or carrying capacity) for a long time, and when extinction happens, the process first leaves a neighbourbood around that level. Populations can be modeled by for example branching processes (such as those treated in \cite{Ham} and \cite{Hog3}). We will consider the Ricker model as an example at the end of the paper; it has been studied in \cite{Hog4}.

In this paper, we will focus on a process that stays in the neighbourhood of the origin for a long time, and then exits from it. Let $\{X_n^\varepsilon\}_{n=0}^\infty$ be a stochastic process defined by the recursion formula
\begin{equation}\label{recursionequation}
X_{n+1}^\varepsilon = f(X_n^\varepsilon) + \varepsilon\xi_{n+1}, X_0 = 0\in R^d,
\end{equation}
where $f$ is a continuous mapping from $R^d$ to itself with a fixed point at the origin, $\{\xi_n\}_{n=1}^\infty$ is a sequence of i.i.d. random variables and $\varepsilon$ is a small positive parameter. The process $\{X_n^\varepsilon\}_{n=0}^\infty$ is a Markov chain on $R^d$. Under suitable assumptions on $\{\xi_n\}_{n=1}^\infty$ and $f$, the chain is positive recurrent and it makes sense to investigate how long it takes for the process to leave a neighbourhood of the origin.

\bigskip

In \cite{Kle}, Klebaner and Liptser showed how the large deviation principle can be used to get an upper bound on the exit time from a set for a process. As an example in that paper, they considered the autoregressive process, defined as in (\ref{recursionequation}) where $f(x) = ax$, and the innovations follow a normal distribution. The upper bound is sharp, since the corresponding lower bound can be found by other methods (\cite{Rut}). A corresponding multivariate case has been studied in \cite{Jun}, and the results were extended to the ARMA model in \cite{Kos}. Exit times for autoregressive processes with other noise distributions have also been studied, in \cite{Hog2}. Different aspects of the exit time problem for linear autoregressive processes are treated in, for example, \cite{Bas}, \cite{Bau}, \cite{Din} and \cite{Nov2}.

\bigskip

It seems natural to continue from the linear case to piecewise linear functions.  Some examples of this type were studied by And\v el (\cite{And}). We will show that the large deviation principle gives explicit (asymptotic) upper bounds on the exit time also in these cases.

\section{Large deviation tools}

In this section we give a summary of how the large deviation principle (LDP) can be used to get an upper bound of the asymptotics of an exit time from a set for a process. The section is based on work by Klebaner and Liptser in \cite{Kle} and by Jung in \cite{Jun}.

\subsection{LDP for a past-dependent process}

The definition of the LDP  used in \cite{Kle} is as follows (this is according to Varadhan's definition in \cite{Var}, with the addition that the rate of speed $q(\varepsilon)$ is a function of $\varepsilon$ as in \cite{Kle}:

\medskip

 Let $\{P_\varepsilon\}$ be a family of probability measures on the Borel subsets of a complete separable metric space $Z$. The family $\{P_\varepsilon\}$ satisfies the large deviation principle with a rate function $I(\cdot)$ if there is a function $I$ from $Z$ into $[0,\infty]$ that satisfies the following conditions: $0\le I(z)\le\infty\;\forall z\in Z$, $I$ is lower semicontinuous, the set $\{z:I(z)\le l\}$ is a compact set in  $Z$ for all $\;l<\infty$ and
\begin{eqnarray*}
& &\limsup_{\varepsilon\rightarrow 0} q(\varepsilon )\log P_\varepsilon (C) \le -\inf_{z\in C}I(z) \;\; \mbox{ for every closed set } C\subset Z \mbox{ and }\\
& &\liminf_{\varepsilon\rightarrow 0} q(\varepsilon )\log P_\varepsilon (G) \ge -\inf_{z\in G}I(z) \;\; \mbox{ for every open set } G\subset Z.
\end{eqnarray*}

\noindent
In \cite{Kle}, Klebaner and Liptser considered a family of processes of the type
\begin{equation}\label{generalprocess}
X_{n+1}^\varepsilon = g(X_n^\varepsilon, \ldots , X_{n-m+1}^\varepsilon , \varepsilon\xi_{n+1}),
\end{equation}
where $g$ is a continuous function on $R$ and $\{\xi_n\}_{n=m}^\infty$ is a sequence of i.i.d. random variables. They gave conditions under which the LDP holds for the family $\varepsilon\xi$ (where $\xi$ is a copy of $\xi_m$) and proved that when $\varepsilon\xi$ obeys an LDP with rate function $I_{\varepsilon\xi}(z)$, it follows that $(X_n^\varepsilon)$ obeys an LDP with a rate function that can be written explicitly using $I_{\varepsilon\xi}(z)$. 

A corresponding result in the multivariate case can be found in \cite{Jun}. Assuming that the function $g$ in (\ref{generalprocess}) is continuous on $R^d$, and that the family $\varepsilon\xi$ obeys an LDP with rate function $I_{\varepsilon\xi}(z)$, the resulting rate function for the family of processes in (\ref{generalprocess}) is
\begin{equation}
I(y_0,y_1,y_2,\ldots ) = \inf_{\genfrac{}{}{0pt}{1}{\genfrac{}{}{0pt}{1}{z_t\in R^d\;\forall t\ge n}{y_t = g(y_{t-1},\ldots , y_{t-n},z_t), t\ge n}}{y_0 = x_0,\ldots , y_{n-1}= x_{n-1}}} \sum_{t=n}^\infty I_{\varepsilon\xi}(z_t).
\end{equation}

\subsection{Upper bound for exit time}

Let the exit time $\tau$ be defined as
\begin{equation}
\tau := \min\{t\ge m: X_n^\varepsilon \notin \Omega\}
\end{equation}
for a set $\Omega$. For the expected exit time it holds that 
\begin{equation}
E_{x_0,\ldots ,x_{m-1}}(\tau ) \le \frac{2M}{\inf_{x_0,\ldots ,x_{m-1}\in \Omega}P_{x_0,\ldots ,x_{m-1}}(\tau\le M)}
\end{equation}
for any set of starting points $x_0,\ldots ,x_{m-1}\in \Omega$ and any integer $M\ge m$ (for details, see \cite{Jun}). If the infimum in the denominator is attained for the starting points $x_0^*,\ldots ,x_{m-1}^*\in\Omega$, the inequality above implies that
\begin{equation}\label{inequality}
\limsup_{\varepsilon\rightarrow 0} q(\varepsilon )\log E_{x_0,\ldots ,x_{m-1}}(\tau) \le -\lim_{\varepsilon\rightarrow 0} q(\varepsilon )\log P_{x_0^*,\ldots ,x_{m-1}^*}(\tau\le M),
\end{equation}
if the right hand side limit exists. Since
\[
P_{x_0^*,\ldots ,x_{m-1}^*}(\tau\le M) = P_{x_0^*,\ldots ,x_{m-1}^*}(X_t^\varepsilon \notin \Omega \mbox{ for some } t\in \{m,\ldots , M\}),
\]
the limit on the right hand side in (\ref{inequality}) may be calculated if we have a large deviation principle for the family of probability measures induced by $\{X_t^\varepsilon\}_{t\ge 0}$ and if the function $f$ and the set $\Omega$ are suitable. 

\subsection{LDP and exit times for process with normally distributed noise}

Consider the case
\begin{equation}\label{processdefinition}
X_{n+1}^\varepsilon = f(X_n^\varepsilon) + \varepsilon\xi_{n+1},
\end{equation}
where $f$ is a continuous function on $\mathbb{R}$ and $\{\xi_n\}_{n\ge 1}$ is an i.i.d. sequence of standard normal random variables. Then $I_{\varepsilon\xi}(z) = \frac{z^2}{2} $ and $g(y_{n-1},\ldots , y_{n-m+1},z_n) = f(y_{n-1}) + z_n$. As Klebaner and Liptser show in \cite{Kle}, this family of processes obeys the large deviation principle with
\begin{equation}
I(y_0,y_1,y_2,\ldots ) =  \frac{1}{2}\sum_{n=1}^\infty (y_n-f(y_{n-1}))^2.
\end{equation}
For the exit time
\begin{equation}
\tau = \min\{n\ge 1: |X_n^\varepsilon |\ge 1\},
\end{equation}
\begin{align}
\nonumber & \limsup_{\varepsilon\rightarrow 0}q(\varepsilon)\log E_{x_0}\tau \le -\lim_{\varepsilon\rightarrow 0} q(\varepsilon) \log P_{x_0^*}(\tau\le M)\\
\nonumber & = \inf_{\genfrac{}{}{0pt}{1}{\genfrac{}{}{0pt}{1}{\max_{1\le n\le M} |y_n|\ge 1}{y_0 = x_0^*}}{}} I(y_0,y_1,y_2,\ldots ) = \inf_{\genfrac{}{}{0pt}{1}{\genfrac{}{}{0pt}{1}{\max_{1\le n\le M} |y_n|\ge 1}{y_0 = x_0^*}}{}} \frac{1}{2}\sum_{n=1}^\infty (y_n-f(y_{n-1}))^2\\
\label{thesum}&= \inf_{1\le N\le M}\left( \inf_{\genfrac{}{}{0pt}{1}{\genfrac{}{}{0pt}{1}{|y_N|\ge 1}{y_0 = x_0^*}}{|y_n|<1, n=1,\ldots ,N}}\frac{1}{2}\sum_{n=1}^N (y_n-f(y_{n-1}))^2\right),
\end{align}
where the upper limit in the sum is $N$, because one can choose $y_n = f(y_{n-1})$ for all $n\ge N+1$ and get the same infimum.

%%%%%%%%%%%%%%%%%%%%%%%%%%%%%%%%%%%%%%%%%%%%%%%%%%%

\section{Minimizing the sum}
In section 2.3. we saw that the asymptotics for the exit times of processes of the type defined in (\ref{processdefinition}) with normally distributed noise is determined by the function $f$ through the infimum of the sum of squares in (\ref{thesum}). In this section we study some properties of these sums for particular classes of autoregression functions $f$.

%%%%%%%%%
%% LEMMA ONE SIDE IS ENOUGH
%%%%%%%%%

\begin{lemma}\label{lemmaonesideenough}
Assume that $f$ is increasing, $f(0) = 0$ (so $f(x) \ge 0$ for $x>0$ and $f(x) \le 0$ for $x < 0$). It is then enough to minimize the sum separately over positive values or negative values:
\begin{align}
\label{infimumlefthandside} & \inf_{\genfrac{}{}{0pt}{1}{\genfrac{}{}{0pt}{1}{|y_N|\ge 1}{y_0 = 0}}{}} \frac{1}{2}\sum_{n=1}^N (y_n-f(y_{n-1}))^2 = \\
\nonumber & \min \left( \inf_{\genfrac{}{}{0pt}{1}{\genfrac{}{}{0pt}{1}{y_N = 1, y_0 = 0}{y_n \ge 0, n= 1,\ldots, N-1}}{}}\frac{1}{2}\sum_{n=1}^N (y_n-f(y_{n-1}))^2,\inf_{\genfrac{}{}{0pt}{1}{\genfrac{}{}{0pt}{1}{y_N = -1, y_0 = 0}{y_n \le 0, n= 1,\ldots, N-1}}{}}\frac{1}{2}\sum_{n=1}^N (y_n-f(y_{n-1}))^2\right)
\end{align}
\end{lemma}
Proof: Assume that the infimum on the left hand side above is attained for $\{y_n^*\}_{n=0}^N$, where $y_0^* = 0$ and $y_N^* = 1$, and let
\begin{equation*}
S^* = \sum_{n=1}^N (y_n^*-f(y_{n-1}^*))^2.
\end{equation*}
Then $S^*\le 1$. We will show by induction that $y_n^*\ge 0$ for $n=1,\ldots , N-1$. We show first that $y_{N-1}^*\ge 0$. Make the counter-assumption that $y_{N-1}^*<0$. Then $f(y_{N-1}^*)\le 0$ and $(1-f(y_{N-1}))^2\ge 1$. Also,
\begin{equation*}
S^* \ge (1-f(y_{N-1}))^2 + (y_L^*)^2,
\end{equation*}
where $L = \min\{i| y_i \neq 0\}$. It follows that $S^* > 1$, which is a contradiction. Thus, the counter-assumption is false and $y_{N-1}^*\ge 0$.

Now, assume that $y_N^*, y_{N-1}^*,\ldots , y_{N-K+1}^*\ge 0$ for some $K<N$. Make the counter-assumption that $y_{N-K}^*<0$. Then $f(y_{N-K}^*)\le 0$ and
\begin{align*}
S^* = & (1-f(y_{N-1}^*))^2 +\ldots + (y_{N-K+2}^*-f(y_{N-K+1}^*))^2\\
& + (y_{N-K+1}^*-f(y_{N-K}^*))^2 + (y_{N-K}^*-f(y_{N-K-1}^*))^2 + \ldots + y_1^2\\
\ge &(1-f(y_{N-1}^*))^2 +\ldots + (y_{N-K+2}^*-f(y_{N-K+1}^*))^2\\
 & + (y_{N-K+1}^*-0))^2 + 0+\ldots + 0.
\end{align*}
In fact, the inequality above is a strict inequality, because
\begin{align*}
S^* \ge & (1-f(y_{N-1}^*))^2 +\ldots + (y_{N-K+2}^*-f(y_{N-K+1}^*))^2\\
 & + (y_{N-K+1}^*-0))^2 + (y_L^*)^2
\end{align*}
where $L = \min\{i\le N-K| y_i^*\neq 0\}$. Thus, $S^*$ is not the minimal sum, which is a contradiction. It follows that $y_{N-K}^*\ge 0$.

If we assume instead that the infimum on the left hand side in (\ref{infimumlefthandside}) is attained for a sequence $\{y_n^\prime\}_{n=0}^N$, where $y_0^\prime = 0$ and $y_N^\prime = -1$, one can show that $y_n^\prime\le 0$ for $n=1,\ldots , N$ in a similar way.

%%%%%%%%%%%%%%
% LEMMA ONLY POSITIVE VALUES FOR ODD FUNCTIONS
%%%%%%%%%%%%%
\begin{lemma}
If $f$ is increasing, $f(0) = 0$ and $f$ is an odd function so that $f(-x) = -f(x)$, we can minimize over only positive values:
\begin{equation}
 \inf_{\genfrac{}{}{0pt}{1}{\genfrac{}{}{0pt}{1}{|y_N|\ge 1}{y_0 = 0}}{}} \frac{1}{2}\sum_{n=1}^N (y_n-f(y_{n-1}))^2 = 
 \inf_{\genfrac{}{}{0pt}{1}{\genfrac{}{}{0pt}{1}{y_N = 1, y_0 = 0}{y_n \ge 0, n= 1,\ldots, N-1}}{}}\frac{1}{2}\sum_{n=1}^N (y_n-f(y_{n-1}))^2
\end{equation}
\end{lemma}
Proof: This follows immediately from lemma \ref{lemmaonesideenough}, since the two infima on the right hand side in (\ref{infimumlefthandside}) have the same value. (One could just as well minimize over only negative values.)

%% LEMMA INCREASING SEQUENCES

\begin{lemma}\label{lemmaincreasing}
If $f$ is increasing, $f(0) = 0$ and $|f(x)|<|x|$ on $(-1,1)\backslash\{0\}$, it is optimal to minimize over positive values and increasing sequences or negative values and decreasing sequences:
\begin{equation}\label{incr_sec}
\inf_{\genfrac{}{}{0pt}{1}{\genfrac{}{}{0pt}{1}{y_N = 1, y_0 = 0}{y_n \ge 0, n= 1,\ldots, N-1}}{}}\frac{1}{2}\sum_{n=1}^N (y_n-f(y_{n-1}))^2 = \inf_{\genfrac{}{}{0pt}{1}{\genfrac{}{}{0pt}{1}{y_N = 1, y_0 = 0}{y_0\le y_1 \le y_2\le \ldots \le y_N}}{}}\frac{1}{2}\sum_{n=1}^N (y_n-f(y_{n-1}))^2
\end{equation}
and
\begin{equation}\label{decr_sec}
\inf_{\genfrac{}{}{0pt}{1}{\genfrac{}{}{0pt}{1}{y_N = 1, y_0 = 0}{y_n \le 0, n= 1,\ldots, N-1}}{}}\frac{1}{2}\sum_{n=1}^N (y_n-f(y_{n-1}))^2 = \inf_{\genfrac{}{}{0pt}{1}{\genfrac{}{}{0pt}{1}{y_N = -1, y_0 = 0}{y_0\ge y_1 \ge y_2\ge \ldots \ge y_N}}{}}\frac{1}{2}\sum_{n=1}^N (y_n-f(y_{n-1}))^2
\end{equation}
\end{lemma}

Proof: We prove equality (\ref{incr_sec}) ((\ref{decr_sec}) can be proven in a similar manner). Assume that the sum $\sum_{n=1}^N (y_n-f(y_{n-1}))^2$ is minimized by the sequence $\{y_n^*\}_{n=0}^N$, where $y_0^* = 0$, $y_N^* = 1$ and $y_n^*\in [0,1]$ for $n=1,\ldots , N-1$. We want to show that $y_N^*\ge y_{N-1}^* \ge \ldots \ge y_1^* \ge y_0^*$. We use induction. It is given that $y_N^*\ge y_{N-1}^*$. Assume that  $y_N^*\ge y_{N-1}^* \ge \ldots \ge y_{N-k+1}^*$ for some $k$. We show that $y_{N-k+1}^*\ge y_{N-k}^*$.

Now, if $y_{N-k+1}^* = 0$, it is clear that the minimum of the sum is attained for $y_{N-k}^* = y_{N-k-1}^* = \ldots = y_1^* = y_0^* = 0$. Then it is clear that $y_{N-k}^*\le y_{N-k+1}^*$. 

If $y_{N-k+1}^* >0$, make the counter-assumption that $y_{N-k+1}^* < y_{N-k}^*$. Then $y_{N-k}^*\in (y_{N-m}^*, y_{N-m+1}^*]$ for some $m\in\{1,\ldots , k-1\}$. Then we have
\begin{align}
\nonumber \sum_{n=1}^N (y_n^*-f(y_{n-1}^*))^2  = &(1-f(y_{N-1}^*))^2 + \ldots + (y_{N-m+2}^*-f(y_{N-m+1}^*))^2\\
\nonumber & + (y_{N-m+1}^*-f(y_{N-m}^*))^2 + (y_{N-m}^*-f(y_{N-m-1}^*))^2\\
\nonumber & + \ldots + (y_{N-k+1}^*-f(y_{N-k}^*))^2 + (y_{N-k}^*-f(y_{N-k-1}^*))^2\\
\nonumber & + \ldots + (y_{2}^*-f(y_{1}^*))^2 + (y_{1}^*)^2\\
\label{suminequality} \ge & (1-f(y_{N-1}^*))^2 + \ldots + (y_{N-m+2}^*-f(y_{N-m+1}^*))^2\\
\nonumber & + (y_{N-m+1}^*-f(y_{N-k}^*))^2 + 0 + \ldots + 0\\
\nonumber & + (y_{N-k}^*-f(y_{N-k-1}^*))^2 + \ldots + (y_1^*)^2,
\end{align}
because $y_{N-m+1}^*-f(y_{N-m}^*) \ge y_{N-m+1}^*-f(y_{N-k}^*)$. Equality is attained if 
\begin{equation}\label{minimumattainedif}
f(y_{N-m}^*) = f(y_{N-k}^*) \mbox{ and } y_{N-m}^* = f(y_{N-m-1}^*), \ldots , y_{N-k+1}^* = f(y_{N-k}^*).
\end{equation}
If this is true, $f$ is constant on the interval $[y_{N-m}^*, y_{N-k}^*]$, and $f(y_{N-m}^*) = y_{N-k+1}^*$. Also,
\begin{equation*}
f(y_{N-m}^*) = f(f(y_{N-m-1}^*)) = \ldots = f(f(\ldots (f( y_{N-k+1}^*)))) < y_{N-k+1}^*,
\end{equation*}
because $y_{N-k+1}^* > 0$. Thus, \ref{minimumattainedif} does not hold, which implies that equality is not attained in \ref{suminequality}. Then, the sum $\sum_{n=1}^N (y_n-f(y_{n-1}))^2$ is not minimized by the sequence $\{y_n^*\}_{n=0}^N$, which is a contradiction. Thus, $y_{N-k+1}^* \ge y_{N-k}^*$.

%%%%%%
% LEMMA TO COMPARE TWO FUNCTIONS
%%%%%%

\begin{lemma}\label{comparefunctions}
If $f$ and $g$ are as in lemma \ref{lemmaincreasing} and $|f(x)|\le |g(x)|$ on $[-1,1]$, then
\begin{equation}
 \inf_{\genfrac{}{}{0pt}{1}{\genfrac{}{}{0pt}{1}{|y_N|\ge 1}{y_0 = 0}}{}}\frac{1}{2}\sum_{n=1}^N (y_n-g(y_{n-1}))^2 \le  \inf_{\genfrac{}{}{0pt}{1}{\genfrac{}{}{0pt}{1}{|y_N|\ge 1}{y_0 = 0}}{}}\frac{1}{2}\sum_{n=1}^N (y_n-f(y_{n-1}))^2
\end{equation}
\end{lemma}

Proof: By lemma \ref{lemmaincreasing}, the minimum on the right hand side is attained for an increasing sequence $\{y_n\}_{n=0}^N$: 
\begin{equation}
\inf_{\genfrac{}{}{0pt}{1}{\genfrac{}{}{0pt}{1}{y_N = 1, y_0 = 0}{y_n \ge 0, n= 1,\ldots, N-1}}{}}\frac{1}{2}\sum_{n=1}^N (y_n-f(y_{n-1}))^2 = \inf_{\genfrac{}{}{0pt}{1}{\genfrac{}{}{0pt}{1}{y_N = 1, y_0 = 0}{y_0\le y_1\le \ldots \le y_N}}{}}\frac{1}{2}\sum_{n=1}^N (y_n-f(y_{n-1}))^2.
\end{equation}
Now, for each $n = 1, \ldots, N$,
\begin{equation}
y_{n}-f(y_{n-1}) \ge y_{n}-y_{n-1} \ge 0,
\end{equation}
and the same is true when $f$ is replaced by $g$.
Since $f(x)\le g(x)$ on $[0,1]$, 
\begin{equation}
y_{n}- f(y_{n-1}) \ge y_{n}-g(y_{n-1}),
\end{equation}
and it follows that
\begin{equation}
(y_{n}- f(y_{n-1}))^2 \ge (y_{n}-g(y_{n-1}))^2.
\end{equation}
The statement in the lemma then follows.

%%%%%%%%%%%%%%%%%%%%%%%%%%%%%%%%%%%%%%%%%%
\subsection{Piecewise linear functions, normal distribution}
Consider the process
\begin{equation}
X_{n+1}^\varepsilon = f(X_n^\varepsilon) + \varepsilon\xi_{n+1},
\end{equation}
where $X_0 = 0$, $\{\xi_n\}_{n\ge 1}$ is an i.i.d. sequence of standard normal random variables and $\varepsilon$ is a small positive parameter. We study different types of piecewise linear functions and the results are given as the following examples.

\begin{example}\label{ex34}
Let $f$ be the function
\begin{equation*}
f(x) = \left\{ \begin{array}{cc} 
				a(x+b) & -1\le x\le -b\\
				0 & -b< x <b\\
				a(x-b) & b\le x \le 1,
\end{array}\right.
\end{equation*}

% \begin{equation}
% f(x) = \left\{ \begin{array}{cc} 
% 				a(x+b) & -1\le x\le -b\\
% 				0 & -b< x <b\\
% 				a(x-b) & b\le x \le 1,
% \end{array}\right.
% \end{equation}
where $|a|<1$ and $0\le b < 1$. Then
\begin{equation}
 \limsup_{\varepsilon\rightarrow 0}\varepsilon^2\log E_{x_0}\tau \le \frac{1}{2}\inf_{1\le N\le M}\left( \frac{(1+\frac{a-a^N}{1-a}b)^2}{\frac{1-a^{2N}}{1-a^2}} \right),
\end{equation}
for any $M\ge 1$. 
\end{example}

% \begin{figure}[!ht]
% \centering
% \includegraphics[width = 0.4\textwidth]{bild1.pdf}
% \caption{Illustration of the function in Example \ref{ex34}.}
% \end{figure}

Proof:
We have
\begin{equation}
\inf _{\genfrac{}{}{0pt}{1}{|y_N| \ge 1}{y_0 = 0}}\sum_{n=1}^N (y_n-f(y_{n-1}))^2 = \inf _{\genfrac{}{}{0pt}{1}{|y_N|= 1}{|y_0| = b}}\sum_{n=1}^N (y_n-f(y_{n-1}))^2
%= \frac{(1+\frac{a-a^N}{1-a}b)^2}{\frac{1-a^{2N}}{1-a^2}}
\end{equation}
Let $c_n = a^{N-n}$ for $n=1\ldots N$ and consider the telescoping sum
\begin{equation}
\sum_{n=1}^N c_n (y_n-f(y_{n-1})) = y_N + ab(c_1+\ldots c_N).
\end{equation}
Cauchy-Schwartz inequality gives
\begin{equation}
\left(\sum_{n=1}^N c_n(y_n-f(y_{n-1}))\right)^2\le \left(\sum_{n=1}^N c_n^2\right)\cdot \sum_{n=1}^N (y_n-f(y_{n-1}))^2.
\end{equation} 
Thus
\begin{equation}\label{lowerboundbeforesumming}
\sum_{n=1}^N (y_n-f(y_{n-1}))^2 \ge \frac{(y_N + ab(c_1+\ldots c_N))^2}{\left(\sum_{n=1}^N c_n^2\right)},
\end{equation}
where equality can be attained.
It follows that
\begin{equation}
 \limsup_{\varepsilon\rightarrow 0}\varepsilon^2\log E_{x_0}\tau \le \frac{1}{2}\inf_{1\le N\le M}\left( \frac{(1+\frac{a-a^N}{1-a}b)^2}{\frac{1-a^{2N}}{1-a^2}} \right),
\end{equation}
for any $M\ge 1$. The value of the infimum, as well as for which $N$ it is attained, depends on the choices of $a$ and $b$, and in some cases, the infimum is one. Some numerical examples are given in the following table.

\begin{center}
\begin{tabular}{|c|c|c|c|} 
 \hline
 a & b & Infimum attained for & Value of infimum \\ 
 \hline
 0 & 0.2 & N=1 & 1 \\
 0.4 & 0.2 & N=1 & 1  \\ 
 0.5 & 0.2 & N=2 & 0.968\\
 0.8 & 0.2 & N=3 & 0.8094\\
 0.95 & 0.2 & N= 4 & 0.6888\\
 \hline
\end{tabular}
\end{center}

If we let $0 < b < 1$ and $a=1$ the formula (\ref{lowerboundbeforesumming}) takes the form 
\begin{equation*}
\frac{1}{2} \inf_{1 \le N \le M}  \frac{(1 + (N-1)b)^2}{N},
\end{equation*}
where the optimal $N$ is either $\left \lfloor{\frac{1}{b}}\right \rfloor  - 1$ or $\left \lfloor{\frac{1}{b}}\right \rfloor$.

\begin{example}\label{ex35}
Let $c\in (0,1]$ and let
\begin{equation*}
f(x) = \left\{ \begin{array}{cc} 
				-ac & -1 \le x< -c\\
				ax & -c\le x\le c\\
				ac & c < x \le 1.
\end{array}\right.
\end{equation*}
Then
\begin{equation}
\limsup_{\varepsilon\rightarrow 0}\varepsilon^2 \log E\tau \le \left\{ \begin{array} {l} \frac{1}{2}((1-ac)^2 + (1-a^2)c^2), c\le a\\ \frac{1}{2}(1-a^2), c\ge a.\end{array}\right.
\end{equation}
\end{example}

% \begin{figure}[!ht]
% \centering
% \includegraphics[width = 0.4\textwidth]{bild2.pdf}
% \caption{Illustration of the function in Example \ref{ex35}.}
% \end{figure}

Proof: We want to find the infimum of 
\begin{equation}
\sum_{n=1}^N (y_n-f(y_{n-1}))^2
\end{equation}
when $y_0 = 0$, $y_0 \le y_1\le \ldots \le y_{N-1}\le y_N$ and $y_N = 1$. The set of sequences which we minimize over can be split into two parts: Either $y_{N-1}< c$ or the $d$ last elements in the sequence are larger than $c$:  $y_{N-d}\ge c$. If $y_{N-d}\ge c$,
\begin{equation}
\sum_{n=1}^N (y_n-f(y_{n-1}))^2 = d\cdot(1-ac)^2 + \sum_{n=1}^{N-d} (y_n-f(y_{n-1}))^2,
\end{equation}
where the latter part is smallest if $y_{N-d} = c$. We get
\begin{equation}
\inf_{\genfrac{}{}{0pt}{1}{\genfrac{}{}{0pt}{1}{y_N = 1, y_0 = 0, y_{N-d}\ge c}{y_0\le y_1\le \ldots \le y_N}}{}}\sum_{n=1}^N (y_n-f(y_{n-1}))^2 = d\cdot(1-ac)^2 + (1-a^2)c^2(1-a^{2(N-d)})^{-1}.
\end{equation}
On the other hand, if $y_{N-1}< c$,
\begin{align}
\inf_{\genfrac{}{}{0pt}{1}{\genfrac{}{}{0pt}{1}{y_N = 1, y_0 = 0, y_{N-1}< c}{y_0\le y_1\le \ldots \le y_N}}{}}\sum_{n=1}^N (y_n-f(y_{n-1}))^2 &= \inf_{\genfrac{}{}{0pt}{1}{\genfrac{}{}{0pt}{1}{y_N = 1, y_0 = 0, y_{N-1}< c}{y_0\le y_1\le \ldots \le y_N}}{}}\sum_{n=1}^N (y_n-ay_{n-1})^2\\
 &= (1-a^2)(1-a^{2N})^{-1}.
\end{align}
We get that
\begin{align*}
\inf_{\genfrac{}{}{0pt}{1}{\genfrac{}{}{0pt}{1}{y_N = 1, y_0 = 0}{y_0\le y_1\le \ldots \le y_N}}{}} & \sum_{n=1}^N (y_n-f(y_{n-1}))^2 \\
& = \min\{d\cdot(1-ac)^2 + (1-a^2)c^2(1-a^{2(N-d)})^{-1}, (1-a^2)(1-a^{2N})^{-1}\}
\end{align*}
and that
\begin{align*}
\inf_{1\le N\le M} & \inf_{\genfrac{}{}{0pt}{1}{\genfrac{}{}{0pt}{1}{y_N = 1, y_0 = 0}{y_0\le y_1\le \ldots \le y_N}}{}} \sum_{n=1}^N (y_n-f(y_{n-1}))^2\\
&\le \min\{d\cdot(1-ac)^2 + (1-a^2)c^2(1-a^{2(M-d)})^{-1}, (1-a^2)(1-a^{2M})^{-1}\}.
\end{align*}
Since 
\begin{equation}
\limsup_{\varepsilon\rightarrow 0}\varepsilon^2 \log E\tau \le \inf_{1\le N\le M} \inf_{\genfrac{}{}{0pt}{1}{\genfrac{}{}{0pt}{1}{y_N = 1, y_0 = 0}{y_0\le y_1\le \ldots \le y_N}}{}} \sum_{n=1}^N (y_n-f(y_{n-1}))^2
\end{equation}
for any positive integer $M$, we can let $M$ be arbitrarily large. It is also optimal to choose $d=1$. We have
\begin{equation}
\limsup_{\varepsilon\rightarrow 0}\varepsilon^2 \log E\tau \le \left\{ \begin{array} {l} \frac{1}{2}((1-ac)^2 + (1-a^2)c^2), c\le a\\ \frac{1}{2}(1-a^2), c\ge a.\end{array}\right.
\end{equation}

%%%%%%%%%%%%%%%%%%%%%%%%%%%%%%%%%%%%%%%%%%%%
\begin{example}
Let 
\begin{equation*}
f(x) = \left\{ \begin{array}{cc} 
				0 & -1 \le x< 0\\
				-ax & 0\le x\le 1.
\end{array}\right.
\end{equation*}
Then 
\begin{equation*}
\limsup_{\varepsilon\rightarrow 0}\varepsilon^2 \log E\tau \le \frac{1}{2}\frac{1}{1+a^2}.
\end{equation*}
\end{example}
We consider
\begin{equation*}
S_N = \inf_{\genfrac{}{}{0pt}{1}{\genfrac{}{}{0pt}{1}{y_N = 1, y_0 = 0}{}}{}} \sum_{n=1}^N (y_n-f(y_{n-1}))^2.
\end{equation*}
If $N=1$, $S_N = 1$. Let $N\ge 2$. If $y_N = 1$, $S_N \ge 1$. Also, if $y_N = -1$ and $y_{N-1}<0$, $S_N\ge 1$. When $y_N = -1$ and $y_{N-1}\ge 0$, it is optimal to have $y_{N-2} = \ldots = y_1 = y_0 = 0$. The minimum of 
\begin{equation*}
(-1-f(x))^2 + x^2
\end{equation*}
for $x\in[0,1]$ is $\frac{1}{1+a^2}$ (it is attained for $x = \frac{a}{1+a^2}$). Thus, 
\begin{equation*}
\inf_{1\le N\le M} \inf_{\genfrac{}{}{0pt}{1}{\genfrac{}{}{0pt}{1}{y_N = 1, y_0 = 0}{y_0\le y_1\le \ldots \le y_N}}{}} \sum_{n=1}^N (y_n-f(y_{n-1}))^2 = \frac{1}{1+a^2}.
\end{equation*}

\begin{example}
Let 
\begin{equation*}
f(x) = \left\{ \begin{array}{cc} 
				-bx & -1 \le x< 0\\
				-ax & 0\le x\le 1,
\end{array}\right.
\end{equation*}
where $0<a<1$, $0<b<1$. Then
\begin{equation*}
\limsup_{\varepsilon\rightarrow 0}\varepsilon^2 \log E\tau \le \frac{1}{2} \min\left( \frac{1-(ab)^2}{1+a^2},\frac{1-(ab)^2}{1+b^2}\right).
\end{equation*}
\end{example}
Proof: Let 
\begin{equation*}
S = \sum_{n=1}^N (y_n-f(y_{n-1})^2.
\end{equation*}
We consider the cases $y_N = 1$ and $y_N= -1$ separately. First, let $y_N= -1$. Clearly, the sum $S$ is smallest if the sequence $\{y_n\}_{n=1}^N$ has alternating signs: $y_{N-1} > 0$, $y_{N-2} <0, \ldots$ Differentiating $S$ with respect to $y_n$ gives
\begin{equation*}
(y_n-f(y_{n-1})) = f^\prime (y_n) (y_{n+1}-f(y_n)),
\end{equation*}
where $f^\prime (y_n) = -b$ if $y_n<0$ and $f^\prime (y_n) = -a$ if $y_n>0$. Let $s_n := y_n-f(y_{n-1})$. Then, if $N= 2M$ so that $N$ is even, it is optimal to have $y_1>0$. Then
\begin{align*}
s_{2k} = -a^{-k}b^{-k+1}s_1\\
s_{2k+1} = a^{-k}b^{-k}s_1,
\end{align*}
and it follows that 
\begin{equation*}
y_N = -s_1[a^Mb^{M-1} + a^{M-2}b^{M-1} + a^{M-2}b^{M-3}+ a^{M-4}b^{M-3}+\ldots + a^{-M+2}b^{-M+1}+ a^{-M}b^{-M+1}]
\end{equation*}
This is a sum of two geometric sums. Since $y_N = -1$, we get
\begin{equation*}
s_1 = \frac{1}{(1+a^2)(ab)^Na^{-2}b^{-1}}\cdot \frac{(ab)^{-2}-1}{(ab)^{-2N}-1}.
\end{equation*}
The value of the sum is
\begin{equation*}
S = \sum_{n=1}^N s_n^2 = \frac{1-(ab)^2}{1+a^2}\cdot \frac{1}{1-(ab)^{2N}}. 
\end{equation*}
If $N$ is odd, it is optimal to have $y_1<0$. Then the same calculations as above follow, if $a$ is replaced by $b$ and vice versa. It follows that the value of the sum is 
\begin{equation*}
S = \sum_{n=1}^N s_n^2 = \frac{1-(ab)^2}{1+b^2}\cdot \frac{1}{1-(ab)^{2N}}. 
\end{equation*}
The cases when $y_N = 1$ and $N$ is odd or even, give the same values of the sum. The smallest values are attained when $N$ is large, and the statement in the example follows.

\bigskip

We note that in the case $a=b$, the result coincides with the earlier result for the autoregressive process, as expected.

\begin{example}\label{absolutevalueexample}
Let 
\begin{equation*}
f(x) = a |x|, -1\le x\le 1
\end{equation*}
for $|a|<1$.
Then 
\begin{equation}\label{upperboundabsolutevalue}
\limsup_{\varepsilon\rightarrow 0}\varepsilon^2 \log E\tau \le \frac{1-a^2}{2}.
\end{equation}
\end{example}
Proof: If $0\le a < 1$,
\begin{equation}
 \inf_{\genfrac{}{}{0pt}{1}{\genfrac{}{}{0pt}{1}{|y_N|\ge 1}{y_0 = 0}}{}} \frac{1}{2}\sum_{n=1}^N (y_n-a|y_{n-1}|)^2 = 
 \inf_{\genfrac{}{}{0pt}{1}{\genfrac{}{}{0pt}{1}{y_N = 1, y_0 = 0}{y_n \ge 0, n= 1,\ldots, N-1}}{}}\frac{1}{2}\sum_{n=1}^N (y_n-ay_{n-1})^2
\end{equation}
and we have the same infimum as in the autoregressive case (see \cite{Kle} for details). If $-1< a \le 0$, 
\begin{equation}
 \inf_{\genfrac{}{}{0pt}{1}{\genfrac{}{}{0pt}{1}{|y_N|\ge 1}{y_0 = 0}}{}} \frac{1}{2}\sum_{n=1}^N (y_n-a|y_{n-1}|)^2 = 
 \inf_{\genfrac{}{}{0pt}{1}{\genfrac{}{}{0pt}{1}{y_N = -1, y_0 = 0}{y_n \le 0, n= 1,\ldots, N-1}}{}}\frac{1}{2}\sum_{n=1}^N (y_n-ay_{n-1})^2,
\end{equation}
which also gives the same result as in the autoregressive case. Thus, we get (\ref{upperboundabsolutevalue}).
% \section{Poisson distributed noise}
% In \ref{Kle} Klebaner and Liptser also gave as an example the case when $\xi$ is a  
% Poisson random variable with parameter 1 and the cumulant function $H(t) =e^t+e^{-t}-2$.

\section{Quadratic functions}
Consider now the process
\begin{equation}
X_{n+1}^\varepsilon = f(X_n^\varepsilon) + \varepsilon\xi_{n+1},
\end{equation}
where $X_0 = 0$, $\{\xi_n\}_{n\ge 1}$ is an i.i.d. sequence of standard normal random variables, $\varepsilon$ is a small positive parameter and $f(x) = ax^2$. 

\begin{example}
When $f(x) = ax^2$ and $0\le a\le 0.5$,
\begin{equation*}
\limsup_{\varepsilon\rightarrow 0}\varepsilon^2 \log E\tau \le \frac{1}{2}.
\end{equation*}
If $a\ge 0.5$,
\begin{equation}\label{quadraticabigger}
\limsup_{\varepsilon\rightarrow 0}\varepsilon^2 \log E\tau \le \frac{1}{2}\left(\frac{1}{a}-\frac{1}{4a^2}\right).
\end{equation}
\end{example}
Proof: Since $f(x)\ge 0$ and $f$ is even, it is optimal to have $y_N = 1$ and $y_n\ge 0\forall n$. We have
\begin{align*}
\sum_{n=1}^N &(y_n-f(y_n))^2 = \sum_{n=1}^N (y_n-ay_{n-1}^2)^2\\
&= 1+y_{n-1}^2(1-2a) + y_{n-2}^2(1-2ay_{n-1}) + \ldots + y_1^2(1-2ay_2)\\
& + a^2(y_{n-1}^4 + y_{n-2}^4 + \ldots + y_2^2) \ge 1
\end{align*}
when $0\le a\le 0.5$, and equality is achieved by putting $y_1 = y_2 = \ldots = y_{N-1} = 0$.

For $a\ge 0.5$, 
\begin{equation*}
\inf_{1\le N\le M}\left(\inf_{\genfrac{}{}{0pt}{1}{\genfrac{}{}{0pt}{1}{|y_N| \ge 1, y_0 = 0}{|y_n|< 1, n=1,2,\ldots , n-1}}{}}\frac{1}{2}\sum_{n=1}^N (y_n-ay_{n-1}^2)^2 \right) \le \inf_{\genfrac{}{}{0pt}{1}{\genfrac{}{}{0pt}{1}{y_2 = 1, y_0 = 0}{|y_1|< 1, }}{}}\frac{1}{2}\sum_{n=1}^2 (y_n-ay_{n-1}^2)^2,
\end{equation*}
where the infimum on right hand side is  $\frac{1}{2}(\frac{1}{a}-\frac{1}{4a^2})$ (it is attained for $y_1 = \sqrt{\frac{1}{a}-\frac{1}{2a^2}}$). This gives the upper bound in (\ref{quadraticabigger}). This is not necessarily the best upper bound for all $a\ge 0.5$, since we have only calculated the infimum for $N=2$.

%%%%%%%%%%%%%%%%%%%%%%%%%%%%%%%%%%%%%%%%%%%%%%%%%%
\section{Connection to stationary distribution}
Mark Kac proved in 1947 that the mean return time of a discrete Markov
chain to a point $x$ is the reciprocal of the invariant probability $\pi(x)$. 

In \cite{Hog}, we explored this idea by comparing the exit time for the process defined by the stochastic difference equation
\begin{equation}\label{stochdiffequation}
X_{n+1}^\varepsilon = f(X_n^\varepsilon) + \varepsilon\xi_{n+1}
\end{equation}
and the return time to a certain set for the same process. We get an upper bound for the asymptotics of the exit time, and that bound is the reciprocal of the stationary distribution of the process, evaluated at the point where the level curve of the stationary distribution touches the boundary of the set (or interval in the univariate case) from which the process exits.

In the univariate case with $f(x) = ax$ (that is, for the autoregressive process), this method gives the same upper bound as the LDP method:
\begin{equation*}
\limsup\varepsilon^2\log E\tau \le \frac{1-a^2}{2}.
\end{equation*}
These methods only give an upper bound, but we know that the bound is sharp in the autoregressive case; the corresponding lower bound can be shown by other methods (\cite{Rut}, where Novikov's martingale method (\cite{Nov} and \cite{Nov2}) was used, and \cite{Jun}).

In section 3, we saw several examples of processes of the type in (\ref{stochdiffequation}) where $f$ was a piecewise linear function. For these examples, it is not straightforward to derive expressions for the stationary distributions of the processes, but in one of the cases, we can see the connection to the stationary distribution explicitly:

And\v el (\cite{And}) considered the process
\begin{equation*}
X_{n+1}^\varepsilon = f(X_n^\varepsilon) + \varepsilon\xi_{n+1},
\end{equation*}
where $f(x) = -|ax|$ for $-1\le x\le 1$,
and found the stationary distribution
\begin{equation*}
\frac{1}{\varepsilon}(\frac{2(1-a^2)}{\pi})^{1/2} \exp(-\frac{1}{2}(1-a^2)\frac{x^2}{\varepsilon^2}) \Phi (\frac{-ax}{\varepsilon}).
\end{equation*}

It holds that
\begin{equation*}
-\varepsilon^2\log (\frac{1}{\varepsilon}(\frac{2(1-a^2)}{\pi})^{1/2} \exp(-\frac{1}{2}(1-a^2)\frac{x^2}{\varepsilon^2}) \Phi (\frac{-ax}{\varepsilon})) \rightarrow \frac{1-a^2}{2},
\end{equation*}
at the point $x=-1$. This limit is the upper bound found in Example \ref{absolutevalueexample}. 

The upper bound in Example \ref{absolutevalueexample} is sharp; one could use Novikov's martingale method for this process to get the corresponding lower bound (the calculations are the same as for the autoregressive process, see \cite{Rut}).
%%%%%%%%%%%%%%%%%%%%%%%%%%%%%%%%%%%%%%%%%%%%%%%%%
\section{The Ricker model}

The classical deterministic Ricker model is defined by
\begin{equation*}
x_{n+1} = x_ne^{r-\gamma x_n}, n= 0,1,2,\ldots,
\end{equation*}
where $x_t$ represents the size or density of a population at time $t$, $r >0$ models the growth rate and $\gamma > 0$ is an environmental factor (\cite{Hog4}). By suitably renorming the population, we may take $\gamma = 1$. The model then has a fixed point at $x=r$ (and one at $x=0$). If we introduce stochasticity in the model by adding white noise, and move the fixed point $x=r$ to the origin, we have the process
\begin{equation*}
X_{n+1}^\varepsilon = f(X_n^\varepsilon) + \varepsilon\xi_{n+1},
\end{equation*}
where $f(x) = (x+r)e^{-x}-r$. We can examine the time until the process exits from a suitable neighbourhood of the origin. 

\begin{figure}[h]%
    \centering
    \subfloat[]{{\includegraphics[width=5cm]{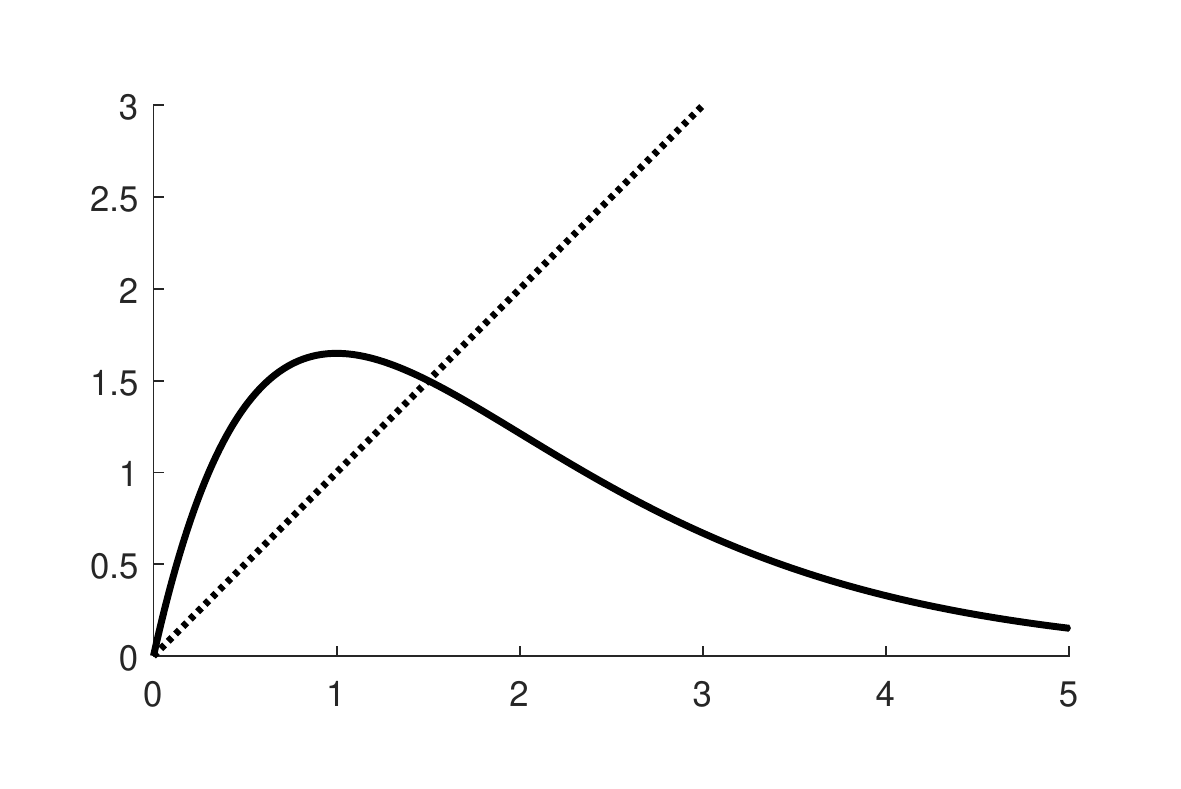} }}%
    \qquad
    \subfloat[]{{\includegraphics[width=5cm]{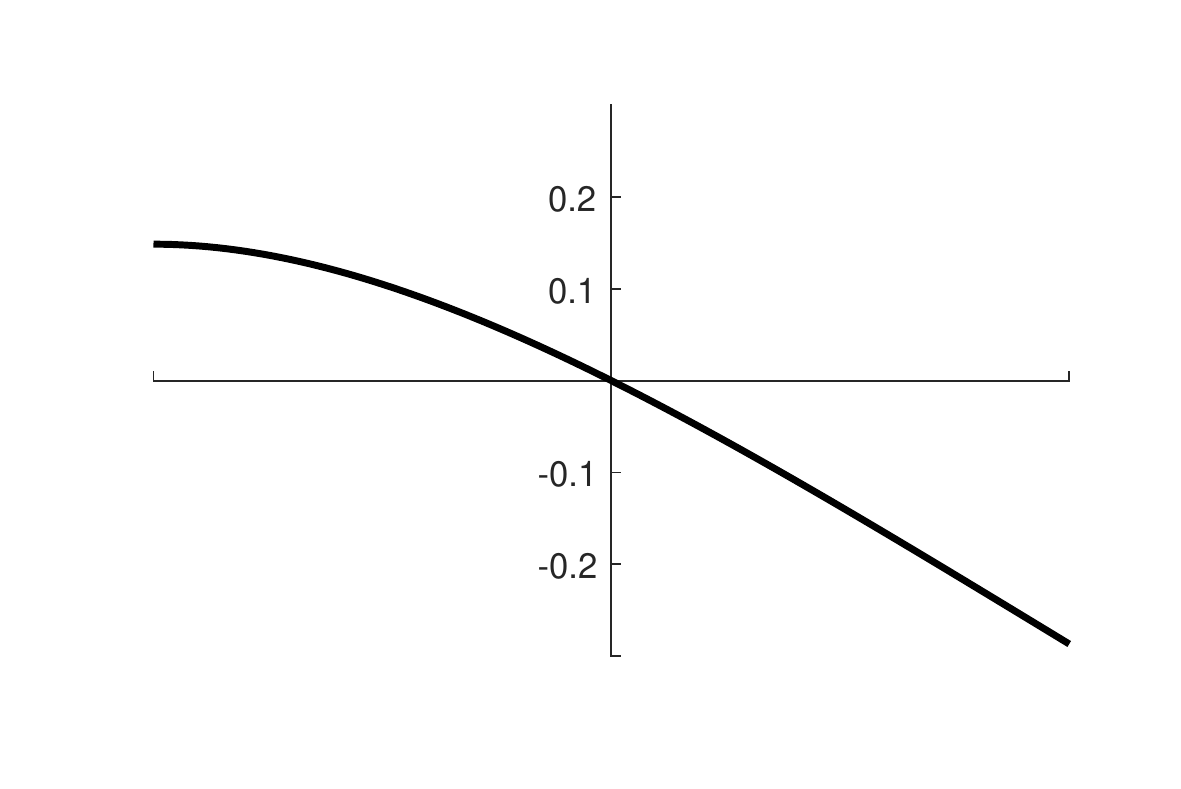} }}%
    \caption{The function $f(x) = x e^{r-x}$ for $r = 1.5$. We have a fixed point at $x= r$. On the right, we see a part of the plot with the fixed point moved to the origin.}%
    \label{fig:example}%
\end{figure}

If $r = 1.5$, consider for example the time until exit from the interval $[-0.5, 0.5]$. Numerical calculations of the infimum
\begin{equation*}
\inf_{1\le N\le M}\left( \inf_{\genfrac{}{}{0pt}{1}{\genfrac{}{}{0pt}{1}{|y_N|\ge 0.5}{y_0 = 0}}{|y_n|<0.5, n=1,\ldots ,N}}\frac{1}{2}\sum_{n=1}^N (y_n-f(y_{n-1}))^2\right)
\end{equation*}
give the approximate value 0.09, so that 
\begin{equation*}
\limsup_{\varepsilon\rightarrow 0}\varepsilon^2\log E\tau \lessapprox 0.09.
\end{equation*}
The derivative of $f$ at the origin is $1-r$, so a suitable linear approximation of the function $f$ is $l(x) = -0.5 x$. By replacing $f$ by $l$, we have an autoregressive process and the upper bound
\begin{equation*}
\limsup_{\varepsilon\rightarrow 0}\varepsilon^2\log E\tau \le 0.5^2\frac{1-0.5^2}{2} = 0.09375.
\end{equation*}
If $r = 0.6$, a suitable interval to consider is $[-0.4,0.4]$. Numerical calculations of the infimum
\begin{equation*}
\inf_{1\le N\le M}\left( \inf_{\genfrac{}{}{0pt}{1}{\genfrac{}{}{0pt}{1}{|y_N|\ge 0.4}{y_0 = 0}}{|y_n|<0.4, n=1,\ldots ,N}}\frac{1}{2}\sum_{n=1}^N (y_n-f(y_{n-1}))^2\right)
\end{equation*}
give the approximate value and upper bound 0.055. The corresponding value of the upper bound if we first make a linear approximation of the function is 
\begin{equation*}
0.4^2\frac{1-0.4^2}{2} = 0.0672.
\end{equation*}

\begin{figure}[ht]%
    \centering
    \subfloat[]{{\includegraphics[width=5cm]{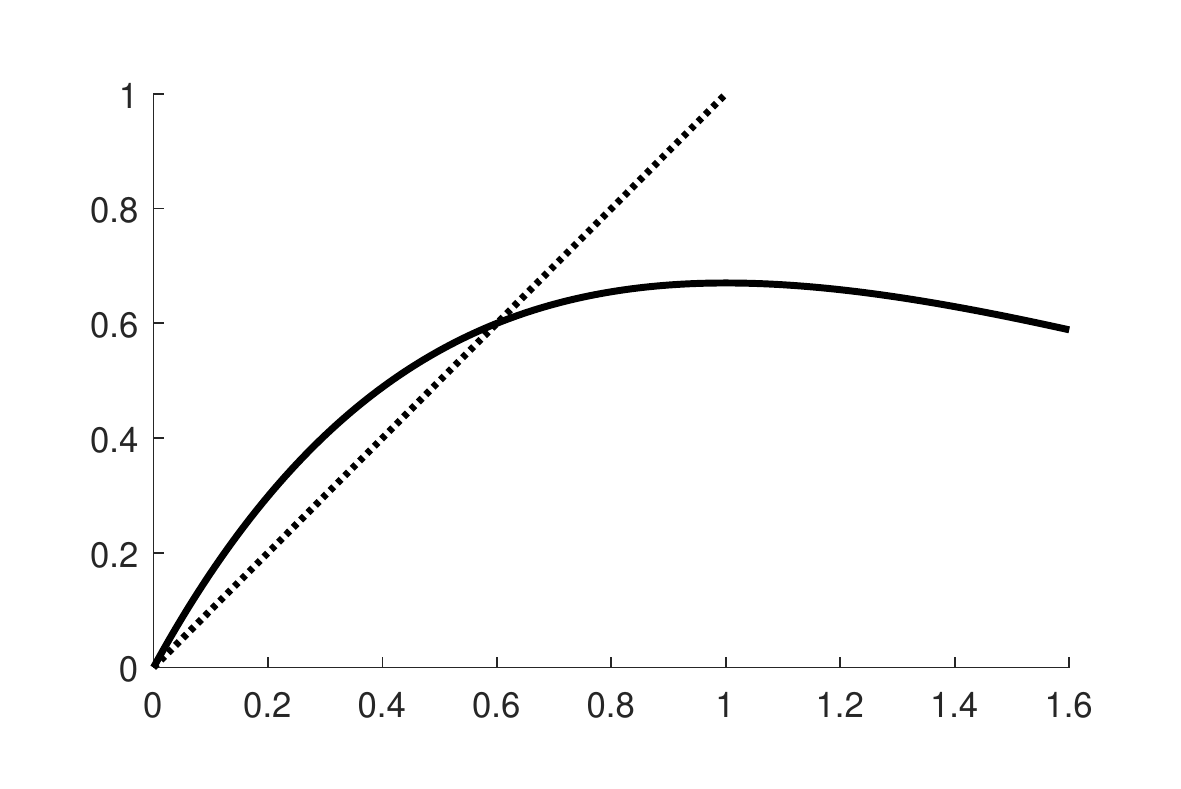} }}%
    \qquad
    \subfloat[]{{\includegraphics[width=5cm]{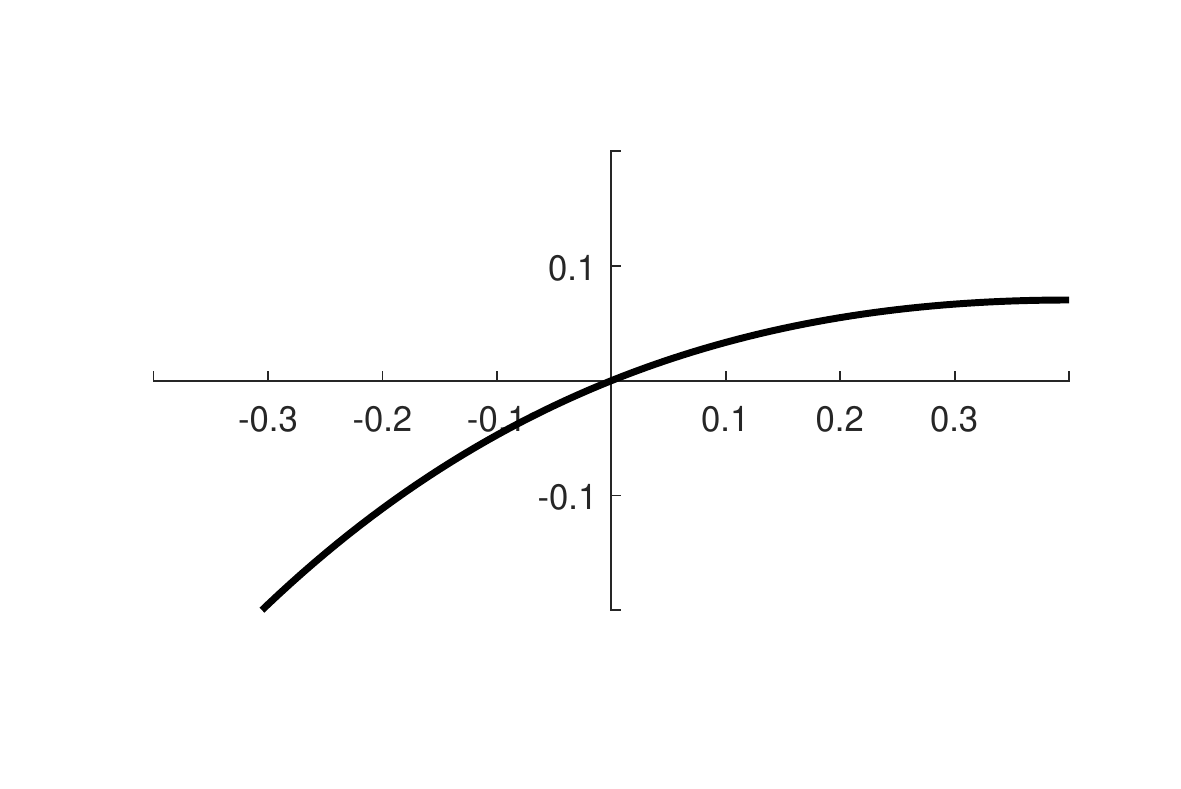} }}%
    \caption{The function $f(x) = x e^{r-x}$ for $r = 0.6$. We have a fixed point at $x= r$. On the right, we see a part of the plot with the fixed point moved to the origin.}%
%     \label{fig:example}%
\end{figure}
We note that a linear approximation might give good enough approximations of the upper bounds, in cases when the neighbourhood around the origin is chosen to be rather small.

\section{Other innovation distributions}

Other types of distributions have been studied before, in \cite{Hog2}. There, other methods were used, not the large deviation principle. The following two cases were treated: If $\xi$ has a Laplace(b)-distribution,
\begin{equation*}
\lim_{\varepsilon\rightarrow 0}\varepsilon\log E\tau  = \frac{1}{b},
\end{equation*}
where we actually have equality and not only an upper bound. If $\xi$ has a Cauchy(0,1)-distribution, then
\begin{equation*}
\limsup_{\varepsilon\rightarrow 0}\varepsilon E\tau  \le \frac{\pi}{2}.
\end{equation*}
(The proofs can be found in \cite{Hog2}.)

The starting point of our investigations concerning autoregressive processes with Gaussian innovations was the example Klebaner and Liptser used in \cite{Kle}. In that paper, another example was also mentioned: The case when $\xi = \xi_1-\xi_2$, a difference between two Poisson($\lambda$)-distributions. Klebaner and Liptser showed that the LDP holds in that case with 
\begin{equation*}
q(\varepsilon) = \frac{\varepsilon}{|\log(\varepsilon)|} \mbox{ and } I(v) = |v|.
\end{equation*}
If we consider a process of the type
\begin{equation}
X_{n+1}^\varepsilon = f(X_n^\varepsilon) + \varepsilon\xi_{n+1}, X_0 = 0,
\end{equation}
for this choice of distribution for $\xi$, the rate function for this family of processes is
\begin{equation}
I(y_0,y_1,y_2,\ldots ) =  \frac{1}{2}\sum_{n=1}^\infty |y_n-f(y_{n-1})|.
\end{equation}
This implies that
\begin{equation*}
\limsup_{\varepsilon\rightarrow 0}\frac{\varepsilon}{|\log \varepsilon|}\log E_{x_0}\tau^\varepsilon \le  \inf_{1\le N\le M}\left( \inf_{\genfrac{}{}{0pt}{1}{\genfrac{}{}{0pt}{1}{|y_N|\ge 1}{y_0 = x_0^*}}{}}\lambda\sum_{n=1}^N |y_n-f(y_{n-1})|\right).
\end{equation*}
If $f$ is as in lemma \ref{lemmaincreasing}, the infimum is attained for $N=1$ and the right hand side is then $\lambda$.

\end{document}